\documentclass[10pt]{article}
\usepackage{mathrsfs}
\usepackage{amsmath,amscd,amsthm,amsfonts}
\usepackage[dvips]{graphicx}
\usepackage[dvips]{psfrag}
\oddsidemargin=0.5cm\topmargin=-1cm \numberwithin{equation}{section}

\theoremstyle{plain}
\newtheorem{theorem}{theorem}[section]

\newtheorem{lemma}[theorem]{Lemma}
\newtheorem{sublemma}[theorem]{Sublemma}

\newtheorem{remark}[theorem]{Remark}

\newtheorem*{ThmA}{Theorem A}
\newtheorem*{ThmB}{Theorem B}
\newtheorem*{ThmC}{Theorem C}
\newtheorem*{ThmD}{Theorem D}
\newtheorem*{ThmE}{Theorem E}
\newtheorem*{ThmB$'$}{Theorem B$'$}

\def\ve{\varepsilon} 

\def\id{{\rm id}}
\def\mod{\mathop{\hbox{\rm mod}}}

\def\vol{\mathop{\hbox{\rm vol}}}

\def\disp{\displaystyle}

\begin{document}
\vskip 5.1cm
\title{Quasi-Shadowing for Partially Hyperbolic Diffeomorphisms
}
\author {HUYI HU$^1$, YUNHUA ZHOU$^2$ and YUJUN ZHU$^3$   \\
\small {1. \emph{Department of Mathematics, Michigan State University, East Lansing, MI 48824, USA}}\\
\small {\emph{(e-mail: hhu@math.msu.edu)}}\\
\small {2. \emph{College of Mathematics and Statistics, Chongqing University, Chongqing, 401331, China}}\\
\small{\emph{(e-mail: zhouyh@cqu.edu.cn)}}\\
\small {3. \emph{College of Mathematics and Information Science, Hebei Normal University, }}\\
\small {\emph{Shijiazhuang, 050024, China}}\\
\small{\emph{(e-mail: yjzhu@mail.hebtu.edu.cn)}}}
\date{}
\maketitle
\begin{center}
\begin{minipage}{130mm}
{\bf Abstract}:
A partially hyperbolic diffeomorphism $f$ has the quasi-shadowing property
if for any pseudo orbit $\{x_k\}_{k\in \mathbb{Z}}$,
there is a sequence of points $\{y_k\}_{k\in \mathbb{Z}}$ tracing it in which
$y_{k+1}$ is obtained from $f(y_k)$ by a motion $\tau$ along the
center direction.
We show that any partially hyperbolic diffeomorphism has the quasi-shadowing
property, and if $f$ has a $C^1$ center foliation
then we can require $\tau$ to move the points along the center foliation.
As applications, we show that any partially hyperbolic
diffeomorphism is topologically quasi-stable under $C^0$-perturbation.
When $f$ has a uniformly compact $C^1$ center foliation,
we also give partially hyperbolic diffeomorphism versions of some theorems
which hold for uniformly hyperbolic systems, such as the Anosov closing lemma,
the cloud lemma and the spectral decomposition theorem.

\end{minipage}
\end{center}

\section{Introduction}

The goal of this paper is to study the shadowing properties for
partially hyperbolic systems and to use them to study some topological
properties of the systems shared by hyperbolic systems.
For partially hyperbolic diffeomorphisms, a center direction is allowed
in addition to the hyperbolic directions. The presence of the center
direction permits a very rich type of structure in these systems.
For general theory of partially hyperbolic systems, we refer to \cite{Hirsch},
\cite{Pesin}, \cite{Barreira} and \cite{Bonatti}.
On the other hand, there is still hyperbolic structure
in partially hyperbolic systems, and therefore we may see some
phenomena similar to that of hyperbolic systems.

It is well known that an Anosov diffeomorphism has the shadowing property
(see \cite{Bowen} for example).
A sequence of $\xi=\{x_{k}\}_{-\infty}^{+\infty}$
 is said to be a $\delta$-\emph{pseudo orbit} for $f$, if
$$
\sup_{k\in \mathbb{Z}} d(f(x_{k}),\;x_{k+1})\leq \delta.
$$
If for a $\delta$-pseudo orbit $\xi=\{x_{k}\}_{-\infty}^{+\infty}$
there is a point $x\in M$ such that
$$
d(f^k(x),\;x_{k})\leq \varepsilon \;\;\mbox{for all}\;\; k \in
\mathbb{Z},
$$
then we say that the point $x$ ``$\varepsilon$-\emph{shadows}"
(or ``$\varepsilon$-\emph{traces}") the
$\delta$-pseudo orbit $\xi$.  We say that $f$ has the \emph{shadowing
property} if for any $\varepsilon>0$, there exists $\delta>0$ such
that every $\delta$-pseudo orbit is $\varepsilon$-shadowed by some
point.

In this paper, we shall investigate the ``shadowing" property of
partially hyperbolic systems. Let
$f$ be a partially hyperbolic diffeomorphism. We cannot expect that in general
the shadowing property holds for $f$ because of the
existence of the center direction.
We show in Theorem A that for any pseudo orbit $\{x_k\}_{k\in \mathbb{Z}}$,
there is a sequence of points $\{y_k\}_{k\in \mathbb{Z}}$ tracing it
in which $y_{k+1}$ is obtained from $f(y_k)$ by a motion $\tau$
along the center direction. In this case we say that $f$ has the
\emph{quasi-shadowing property}.
Moreover, if the center foliation $\mathcal{W}^c_f$ of $f$ exists and
is of $C^1$, then we can choose the motion $\tau$ that maps points
along the center leaves. This result is
given in Theorem~B. Theorem B$'$ deals with a particular case of $f$ whose center foliation is of one dimensional, the corresponding map $\tau$ in this case can be determined by a flow along the foliation.

In \cite{Hirsch} and \cite{Pesin}, the notion of pseudo orbits
with respect to the plaque of the center foliation is introduced
to investigate the robustness of the center foliation for
normally hyperbolic and partially hyperbolic systems respectively.
Recently, various shadowing properties are investigated by Bonatti and Bohnet \cite{Bonatti1}, Carrasco\cite{Carrasco} and Kryzhevich and Tikhomirov \cite{Tikhomirov}, etc., for partially hyperbolic diffeomorphisms with some assumptions, such as dynamical coherence and uniform compactness, on the center foliation.
In this paper we show that for any partially hyperbolic diffeomorphism $f$,
without any additional assumption, the quasi-shadowing property holds.
Moreover, the method we use is quite different from that in the papers mentioned above. In fact, to obtain our main results we adapt the unified  ``analytic" method in \cite{HZ}, in which the quasi-stability of partially hyperbolic diffeomorphisms are investigated, to our case.

Shadowing property implies some other interesting properties
in the study of hyperbolic systems.
We can obtain some similar results for partially hyperbolic systems
from the quasi-shadowing property.

It is well known that any Anosov diffeomorphism $f$ on $M$ is
\emph{topologically stable} (\cite{Walters}), that is,
for any homeomorphism $g$ $C^0$-close to $f$,
there exists a surjective continuous map $h$ on $M$ such that
$h\circ g=f\circ h$.
Topological stability for hyperbolic systems can be obtained
by shadowing property (\cite{Walters1}, see also \cite{Pilyugin}).
Similarly, as an application of quasi-shadowing property,
we show in Theorem C that any partially hyperbolic diffeomorphism has
\emph{topological quasi-stability} under $C^0$ perturbation,
that is, for any homeomorphism $g$ $C^0$-close to $f$, there exist a
surjective continuous map $h$ from $M$ to itself and a family of locally
defined maps $\{\tau_x:x\in M\}$, which move points
along the center direction such that
$h\circ g(x)=\tau_{f(x)}\circ f\circ h(x)$ for all $x\in M$.
In particular, if center foliation $\mathcal{W}^c_f$ of $f$ exists and
is of $C^1$, then we can choose the motion $\tau$ maps points
along the center foliation.
For more information about the quasi-stabilities for partially hyperbolic
diffeomorphisms under $C^0$ and $C^1$-perturbations and their applications, we refer to \cite{HZ}.

A notable property for a hyperbolic system is the Anosov closing lemma,
which says that if an orbit returns to a small neighborhood of its
initial position, then there is a periodic orbit nearby.
Consequently, for an Anosov diffeomorphism, the closure of the set of
 periodic orbits is equal to its nonwandering set.
It is natural to imagine that for a partially hyperbolic system,
if an orbit returns, then there is a {\it periodic center leaf} nearby.
We prove it using the quasi-shadowing property.
Further, we obtain in Theorem D that periodic center leaves are dense
in the nonwandering set.
If all the center leaves are uniformly compact, then the closure
of the periodic center leaves is equal to the {\it center nonwandering set}
of the map (see Section 1 for the precise meaning).

By the quasi-shadowing property,
we give some versions of the cloud lemma and then the spectral decomposition theorem
when $f$ has uniformly compact $C^1$ center foliation,
which are generalizations of the corresponding results
for the Axiom A systems (see \cite{Bowen}, \cite{Shub}, for example).
More precisely, we show that the center nonwandering set
can be uniquely split into finite disjoint
{\it center topologically transitive} closed subsets,
and each of which can be uniquely split into finite disjoint sets
which are invariant and is {\it center topologically mixing}
under an iteration of $f$.

This paper is organized as follows.
The statements of results are given in Section~1.
In Section~2 we deal with the quasi-shadowing property for the general case, i.e., there is no assumption on the center bundle, in Theorem~A.
Section~3 is for the case that the center foliation exists and is of $C^1$,
and the proofs of Theorem~B and Theorem B$'$ are given there.
The last three sections concern applications of the quasi-shadowing properties.
We study the topological quasi-stability, the denseness of periodic center leaves
and the spectral decomposition in center nonwandering sets
in Section~4, 5 and 6, respectively.

\section{Definition and statement of results}\label{DSN}

Everywhere in this paper, we assume that $M$ is a smooth
$m$-dimensional compact Riemannian manifold. We denote by $\|\cdot\|$
and $d(\cdot,\cdot)$ the norm on $TM$ and the metric
on $M$ induced by the Riemannian metric, respectively.

A diffeomorphism $f: M\to M$ is said to be
\emph{(uniformly) partially hyperbolic} if there exist
numbers $\lambda,\lambda',\mu$ and $\mu'$ with
$0<\lambda<1<\mu$ and $\lambda<\lambda'\leq \mu'<\mu$,
and an invariant decomposition $T_xM=E_x^s\oplus E_x^c\oplus E_x^u$
\ $\forall x\in M$, such that for any $n\ge 0$,
\begin{eqnarray*}
\|d_xf^nv\|  &\!\!\!\!\!  \le C\lambda^n\|v\|\ \   &\text{as} \ v\in E^s(x), \\
C^{-1}(\lambda')^n\|v\| \le \|d_xf^nv\|  &\!\!\!  \le C(\mu')^n\|v\|
               \ &\text{as} \   v\in E^c(x),\\
C^{-1}\mu^n\,\|v\| \leq  \|d_xf^nv\|&\ &\text{as} \  v\in E^u(x)
\end{eqnarray*}
hold for some number $C>0$.  The subspaces
$E_x^s, E_x^c$ and $E_x^u$ are called
\emph{stable, center} and \emph{unstable} subspace, respectively.
Via a change of Riemannian metric we always assume that $C=1$.
Moreover, for simplicity of the notation, we assume that
$\disp \lambda=\frac{1}{\mu}$.

Since $M$ is compact, we can take a constant $\rho_0>0$ such that for
any $x\in M$, the standard exponential mapping
$\exp_{x}:\{v\in T_{x}M:\|v\|<\rho_0\}\to M$ is a $C^{\infty}$
diffeomorphism to the image. Clearly, we have $d(x, \exp_xv)=\|v\|$
for $v\in T_xM$ with $\|v\|<\rho_0$.
For any diffeomorphism $f: M\to M$, we take $\rho=\rho_f\in (0,\rho_0/2)$
such that for any $x,y\in M$ with $d(f^{-1}(x),y)\le \rho$, $v\in T_yM$
with $\|v\|\le \rho$,
$$d(x, f\circ \exp_y v)\le \rho_0/2.
$$
Reduce $\rho$ if necessary such that both sides in equation \eqref{eqn4ThmA}
and \eqref{eqn4ThmB}, in the proof of Theorem A and Theorem B
respectively, are contained in the set $\{v\in T_{x}M:\|v\|<\rho_0\}$.

For a sequence of points $\{x_k\}_{k\in \mathbb{Z}}$
and a sequence of vectors $\{u_k\in E^c_{x_k}\}_{k\in \mathbb{Z}}$
with $\|u_k\|<\rho$ for any $k\in \mathbb{Z}$,
we define a family of smooth maps $\tau^{(1)}_{x_k}=\tau^{(1)}_{x_k}(\cdot, u_k)$
on $B(x_k,\rho)$, $k\in \mathbb{Z}$,  by
\begin{equation}\label{fdef of tau1}
\tau^{(1)}_{x_k}(y)=\exp_{x_k}(u_k+\exp_{x_k}^{-1}y).
\end{equation}

\begin{ThmA}\label{ThmA}
Let $f: M\to M$ be a partially hyperbolic diffeomorphism. Then $f$ has the
quasi-shadowing property in the following sense: for any $\varepsilon\in
(0,\rho)$ there exists $\delta>0$ such that for any $\delta$-pseudo
orbit $\{x_k\}_{k\in \mathbb{Z}}$ of $f$, there exist a sequence of points $\{y_k\}_{k\in \mathbb{Z}}$ and
a sequence of vectors $\{u_k\in E^c_{x_k}\}_{k\in \mathbb{Z}}$
such that
\begin{equation}\label{eqn1thmA}
d(x_k,y_k)<\varepsilon,
\end{equation}
where
\begin{equation}\label{eqn2thmA}
y_k=\tau_{x_k}^{(1)}(f(y_{k-1})).
\end{equation}

Moreover, $\{y_k\}_{k\in \mathbb{Z}}$ and $\{u_k\}_{k\in \mathbb{Z}}$ can be chosen
uniquely so as to satisfy
\begin{equation}\label{eqn3thmA}
y_k\in \exp_{x_k}(E^s_{x_k}+E^u_{x_k}).
\end{equation}
\end{ThmA}

The above theorem does not require any additional condition,
provided that $f$ is a partially hyperbolic diffeomorphism.
Here $\tau^{(1)}_{x_k}$ is a motion in the center direction for
any $k\in \mathbb{Z}$.
If $f$ has $C^1$ center foliation
$\mathcal{W}^c_f$, then we can make $\tau$ to move along the center
foliation. In this case, we denote for any $\varepsilon>0$,
$\Sigma_\varepsilon(x)=\exp_x(H_x(\varepsilon))$, where $H_x(\varepsilon)$ is
the $\varepsilon$-ball in $E_x^s\oplus E_x^u$. Obviously,
$\Sigma_\varepsilon(x)$ is a smooth disk transversal to $E_x^c$ at $x$.
Since the center foliation
$\mathcal{W}^c_f$ is $C^1$, we can conclude that if $y$ is close enough to $x$, then there is a locally defined map $\tau_x^{(2)}$ on
some neighborhood $U(x)$ of $x$ and a  constant $K_1>1$ independent of $x$ such that for any $y\in U(x)$,
\begin{equation}\label{cond1 of tau2}
\tau_x^{(2)}(y) \in \Sigma_\varepsilon(x)\cap \mathcal{W}^c_f(y)
\end{equation}
and
\begin{equation}\label{cond2 of tau2}
d(\tau_x^{(2)}(y),x)<K_1d(y,x).
\end{equation}

\begin{ThmB}\label{theorem2}
Let $f:M \to M$ be a partially hyperbolic diffeomorphism
with $C^1$ center foliation $\mathcal{W}^c_f$. Then $f$ has the
quasi-shadowing property in the following sense:  for any $\varepsilon\in
(0,\rho)$ there exists $\delta>0$ such that for any $\delta$-pseudo
orbit $\{x_k\}_{k\in \mathbb{Z}}$ of $f$, there exists a sequence of points
$\{y_k\}_{k\in \mathbb{Z}}$
such that
\begin{equation}\label{eqn1thmB}
d(x_k,y_k)<\varepsilon,
\end{equation}
where
\begin{equation}\label{eqn2thmB}
y_k=\tau_{x_k}^{(2)}(f(y_{k-1})).
\end{equation}

Moreover, $\{y_k\}_{k\in \mathbb{Z}}$ can be chosen
uniquely so as to satisfy (\ref{eqn3thmA}).

\end{ThmB}

As a particular case, when the center foliation
$\mathcal{W}^c_f$ is $C^1$ and of
dimension one then we can define $\tau$ more directly.
Let $u$ be the vector field consisting of unit vectors in center direction,
i.e., $\|u(x)\|=1$ for any $x\in M$,
and $\varphi^t$ be the flow generated by $u$.
For a sequence of points $\{x_k\}_{k\in \mathbb{Z}}$
and a sequence of real numbers $\{\tilde{\tau}_k\}_{k\in \mathbb{Z}}$,
denote a sequence of smooth maps
$\tau^{(3)}_{x_k}=\tau^{(3)}_{x_k}(\cdot, \tilde\tau_k)$ of $B(x_k,\rho)$
for any $k\in \mathbb{Z}$ given by
$$
\tau^{(3)}_{x_k}(z)=\varphi^{\tilde\tau_k}(z).
$$

\begin{ThmB$'$}\label{ThmB'}
Let $f:M \to M$ be a partially hyperbolic diffeomorphism
with one dimensional $C^1$ center foliation $\mathcal{W}^c_f$.
Then $f$ has the quasi-shadowing property in the following sense:
for any $\varepsilon\in (0,\rho)$ there exists $\delta>0$
such that for any $\delta$-pseudo
orbit $\{x_k\}_{k\in \mathbb{Z}}$ of $f$, there exist a sequence
of points $\{y_k\}_{k\in \mathbb{Z}}$ and a sequence of real numbers
$\{\tilde{\tau}_k\}_{k\in \mathbb{Z}}$ such that
\begin{equation}\label{eqn1thmB'}
d(x_k,y_k)<\varepsilon,
\end{equation}
where
\begin{equation}\label{eqn2thmB'}
y_k=\tau_{x_k}^{(3)}(f(y_{k-1})).
\end{equation}

Moreover, $\{y_k\}_{k\in \mathbb{Z}}$ can be chosen
uniquely so as to satisfy (\ref{eqn3thmA}).

\end{ThmB$'$}


Now we consider the applications of our results.
The first one is about quasi-stability. In \cite{HZ},  topological quasi-stability is given for  partially hyperbolic diffeomorphisms under $C^0$ perturbation. Now we apply the quasi-shadowing property to give another proof.

\begin{ThmC}\label{ThmC}
Assume that $f: M\to M$ is a partially hyperbolic diffeomorphism.
Then $f$ has topological quasi-stability in the sense that there exists
$\varepsilon_0\in (0,\rho)$ satisfying the following conditions: for any
$\varepsilon\in (0,\varepsilon_0)$ there exists $\delta>0$ such that
for any homeomorphism $g$ of $M$ with $d(f,g)<\delta$
there exist a continuous center section $u=\{u_x\in E^c_x: x\in M\}$
and a surjective continuous map $h:M\to M$ such that
\begin{equation}\label{feqnThmC}
h\circ g(x)=\tau^{(1)}_{g(x)}\circ f\circ h(x),\;\;x\in M.
\end{equation}
In addition, $h$ can be chosen uniquely so as to satisfy the
following conditions:
\begin{equation}\label{cond1}
\begin{split}
&d(h,{\id}_M)<\varepsilon,  \\
&\exp_x^{-1}(h(x))\in E_x^s\oplus E_x^u \ \  \text{for} \ x\in M.
\end{split}
\end{equation}

Moreover, if $f$ has $C^1$ center foliation $\mathcal{W}^c_f$, then there exists $h$ as above such that
(\ref{feqnThmC}) holds with $\tau^{(1)}_{g(x)}$ replaced by $\tau^{(2)}_{g(x)}$. Furthermore, if the above $C^1$ center foliation $\mathcal{W}^c_f$ is of one dimensional, then there exist $h$ as above and a continuous function $\tilde{\tau}$ on $M$ such that (\ref{feqnThmC}) holds with $\tau^{(1)}_{g(x)}$ replaced by $\tau^{(3)}_{g(x)}$.
\end{ThmC}

It is well known that for uniformly hyperbolic systems,
the closing lemma holds and therefore
the periodic points are dense in the nonwandering set.
We can get a similar result for partially hyperbolic systems
by using Theorem B.
In this case, periodic center leaves and center nonwandering leaves
play the role as periodic points and nonwandering points, respectively.

A center leaf $W^c(p)$ is said to be a \emph{periodic center leaf} with period
$n\in \mathbb{N}$ if $W^c(p)=W^c(f^n(p))$.  Denote
$$P^c(f)=\{p\in  M: W^c(p) \text{ is a periodic center leaf}\}.
$$

We say that a center leaf $W^c(x)$ is {\em center nonwandering}
if for any neighborhood $U$ of $W^c(x)$ consisting of center leaves,
there is $n\ge 1$ such that $f^nU\cap U\ne\emptyset$.
We denote the {\em center nonwandering set} of $f$ by
$$\Omega^c(f)=\{x\in M: W^c(x) \text{ is center nonwandering}\}.
$$
It is easy to see that $\Omega^c(f)$ is a closed invariant set
and saturated by $W^c$, i.e.,
$W^c(x)\subset \Omega^c(f)$ if $x\in \Omega^c(f)$.
Also we denote by $\Omega(f)$ the nonwandering set of $f$.
Clearly, $\Omega(f)\subset \Omega^c(f)$.

We say that the center foliation is {\em uniformly compact} if
$$\sup\{\vol (W^c(x)): x\in M\}<+\infty,$$
where $\vol(W^c(x))$ is the Riemannian volume restricted to the submanifold
$W^c(x)$ of $M$.
The partially hyperbolic systems with uniformly compact center foliations were
studied in \cite{Bohnet} and \cite{Carrasco}.

It is easy to see that if the center foliation is uniformly compact, then
a center leaf $W^c(x)$ is center nonwandering if and only if
for any $\delta>0$, there is $y\in M$ and $n\in \mathbb{N}$ such that
\begin{equation}\label{fOmegac}
\max\{d_H(W^c(x),W^c(y)), d_H(W^c(x),W^c(f^ny))\}<\delta,
\end{equation}
where $d_H(\cdot, \cdot)$ denotes the Hausdorff distance given by
$d_H(A,B)= \max\limits_{a\in A}\min\limits_{b\in B}d(a,b)$
for closed subsets $A, B\subset M$.

For any set $S\subset M$, denote $W^c(S)=\cup_{x\in S}W^c(x)$.
With the notions $P^c(f)$ and $\Omega^c(f)$, we can get an analogues
of Anosov closing lemma for partially hyperbolic diffeomorphisms
(see Lemma~\ref{L1ThmD} for details).
Based on the results, we have the following theorem.

\begin{ThmD}\label{ThmD}
For any partially hyperbolic diffeomorphism $f: M\to M$ with $C^1$ center
foliation $\mathcal{W}^c_f$, $\Omega(f) \subset \overline{P^c(f)}$.

Moreover, if the center foliation of $f$ is uniformly compact, then
\begin{equation}\label{fThmD}
\overline{P^c(f)}=\Omega^c(f)=W^c\bigl(\Omega(f)\bigr).
\end{equation}
\end{ThmD}

Further, if a partially hyperbolic diffeomorphism has uniformly compact
$C^1$ center foliation, then we have the cloud lemma (see Lemma~\ref{L1ThmE}),
and therefore we can get the spectral decomposition for $\Omega^c(f)$.


The substitutes of the topological transitivity and the topological mixing are
the center topological transitivity and the center topological mixing respectively.
An $f$-invariant set $S$ is said to be {\em center topologically transitive},
if for any two nonempty open sets $U,V$ in $S$, there is $n\in \mathbb{N}$
such that
$$f^n(W^c(U))\cap V\neq \emptyset.$$
$S$ is said to be {\em center topologically mixing},
if for any two nonempty open sets $U,V$ in $S$,
there is $n_0\in \mathbb{N}$ such that
$$f^n(W^c(U))\cap V\neq \emptyset \qquad \forall n\geq n_0.$$

\begin{ThmE}\label{ThmE}
Let $f:M\to M$ be a partially hyperbolic diffeomorphism
with uniformly compact $C^1$ center foliation.
Then $\Omega^c(f)$ is a union of finite pairwise disjoint closed sets
$$\Omega^c(f)=\Omega^c_1\cup\cdots\cup\Omega^c_k.$$
Moreover, for each $i=1,2, \cdots, k$, $\Omega^c_i$ satisfies that

\begin{enumerate}
\item[(a)] $f(\Omega^c_i)=\Omega^c_i$ and $f|_{\Omega^c_i}$
is center topologically transitive;

\item[(b)] $\Omega^c_i=X_{1,i}\cup \cdots \cup X_{n_i,i}$ such that the $X_{j,i}$ are disjoint closed sets, $f(X_{j,i})=X_{j+1,i}$ for $1\leq j\leq n_i-1$, $f(X_{n_i,i})=X_{1,i}$,
and $f^{n_i}|_{X_{j,i}}$ is center topologically mixing.
\end{enumerate}
\end{ThmE}

We call $\Omega^c_i, i=1,2, \cdots, k$, the {\em center basic sets} of $f$.

\begin{remark}
We mention that if $f$ is dynamically coherent and plaque expansive with respect to the center foliation then the similar results in Theorem B hold (\cite{Tikhomirov}). Therefore, if we replace the $C^1$ smoothness of the center foliation by the weaker conditions, we can obtain the same results in Theorem~C, D and E in a similar strategy.
\end{remark}

\section{Quasi-shadowing for the general case}

We prove Theorem A in this section.

Recall that $\|\cdot\|$ is the norm on $TM$. We define the norm
$\|\cdot\|_1$ on $TM$ by $\|w\|_1=\|u\|+\|v\|$ if $w=u+v\in T_xM$
with $u\in E^c_x$ and $v\in E^u_x\oplus E^s_x$. For any sequence $\{x_k\}_{k\in \mathbb{Z}}$, Denote
$$
\mathfrak{X}=\{w=\{w_k\}_{k\in \mathbb{Z}}:w_k\in T_{x_k}M, k\in \mathbb{Z}\},
$$
$$
\mathfrak{X}^c=\{u=\{u_k\}_{k\in \mathbb{Z}}:u_k\in E_{x_k}^c, k\in \mathbb{Z}\}
$$
and
$$
\mathfrak{X}^{us}=\{v=\{v_k\}_{k\in \mathbb{Z}}:v_k\in E_{x_k}^u\oplus E_{x_k}^s, k\in \mathbb{Z}\}.
$$
For any $w=u+v\in \mathfrak{X}$, where $u\in \mathfrak{X}^c$ and $v\in \mathfrak{X}^{us}$, we also define
$$
\|w\|=\sup_{k\in \mathbb{Z}}\|w_k\|
$$
and
$$
\|w\|_1=\|u\|+\|v\|.
$$
By
triangle inequality and the fact that the angles between $E^c$ and
$E^u\oplus E^s$ are uniformly bounded away from zero, we know that there
exists a constant $L$ such that
\begin{equation}\label{fdefL}
\|w\|\le \|w\|_1\le L\|w\|.
\end{equation}

For any $\ve>0$, we denote
$$
\mathfrak{B}(\ve)=\{w\in \mathfrak{X}: \|w\|\le \ve\},\qquad
\mathfrak{B}^{us}(\ve)=\{w\in \mathfrak{X}^{us}: \|w\|\le \ve\},
$$
$$
\mathfrak{B}_1(\ve)=\{w\in \mathfrak{X}: \|w\|_1\le \ve\}.
$$

We denote
$\Pi^{s}_x:T_{x}M\to E^{s}_x$ be the projection onto
$E^{s}_x$ along $E^c_x\oplus E^u_x$. $\Pi^{c}_x$ and $\Pi^{u}_x$ are
defined in a similar way.

\begin{proof}[Proof of Theorem A]
Given a $\delta$-pseudo orbit $\{x_k\}_{k\in \mathbb{Z}}$ of
$f$. To find a sequence of points $\{y_k\}_{k\in \mathbb{Z}}$ and
a sequence of vectors $\{u_k\in E^c_{x_k}\}_{k\in \mathbb{Z}}$
satisfying (\ref{eqn1thmA}), (\ref{eqn2thmA}) and (\ref{eqn3thmA}), we shall
try to solve the equations
\begin{equation}\label{feqn2ThmA}
y_k=\tau_{x_k}^{(1)}(f(y_{k-1})), \quad  k\in \mathbb{Z},
\end{equation}
for unknown $\{y_k\}_{k\in \mathbb{Z}}$ and $\{u_k\in E^c_{x_k}\}_{k\in \mathbb{Z}}$,
where $\tau_{x}^{(1)}$ is defined in (\ref{fdef of tau1}).
Put $v_k=\exp_{x_k}^{-1}y_k, k\in\mathbb{Z}$.
Then the equations (\ref{feqn2ThmA}) can be written as
$$
v_{k}=\exp_{x_{k}} ^{-1}\tau^{(1)}_{x_{k}}( f\circ \exp_{x_{k-1}}v_{k-1}),
\quad  k\in \mathbb{Z}.
$$
By (\ref{fdef of tau1}), the equations are equivalent to
\begin{equation}\label{eqn4ThmA}
v_{k}=u_{k}+\exp_{x_{k}}^{-1}\circ f\circ \exp_{x_{k-1}}v_{k-1},\quad k\in \mathbb{Z}.
\end{equation}

Define an operator $\beta: \mathfrak{B}^{us}(\rho)\to \mathfrak{X}$
and a linear operator $A: \mathfrak{B}^{us}(\rho)\to \mathfrak{X}^{us}$ by
\begin{equation}\label{betaThmA}
(\beta(v))_{k}=\exp_{x_{k}}^{-1}\circ f\circ \exp_{x_{k-1}}v_{k-1},
\end{equation}
and
\begin{equation}\label{AThmA}
(Av)_k=((A^s+A^u)v)_k=(A_{k-1}^s+A_{k-1}^u)v_{k-1},
\end{equation}
where
\begin{equation}
\begin{split}
A_{k-1}^s=\Pi_{x_k}^s\circ
d_0(\exp_{x_k}^{-1}\circ f\circ \exp_{x_{k-1}})\circ \Pi_{x_{k-1}}^s, \\
A_{k-1}^u=\Pi_{x_k}^u\circ
d_0(\exp_{x_k}^{-1}\circ f\circ \exp_{x_{k-1}})\circ \Pi_{x_{k-1}}^u.
\end{split}
\end{equation}

Let $\eta=\beta-A$.  By (\ref{betaThmA}) and (\ref{AThmA}),
(\ref{eqn4ThmA}) is equivalent to
$$
v=u+Av+\eta(v),
$$
or
$$
v-u-Av=\eta(v).
$$

Define a linear operator $P$ from a neighborhood of $0\in
\mathfrak{X}$ to $\mathfrak{X}$ by
\begin{equation}\label{PThmA}
Pw=-u+(\id_{\mathfrak{X}^{us}}-A)v,
\end{equation}
and then define an operator $\Phi$ from a neighborhood of $0\in
\mathfrak{X}$ to $\mathfrak{X}$ by
$$
\Phi(w)=P^{-1}\eta(v)
$$
for $w=u+v$ in the neighborhood of $0\in \mathfrak{X}$, where
$u\in \mathfrak{X}^c$ and $v\in \mathfrak{X}^{us}$.

Hence, the equations (\ref{eqn4ThmA}) are equivalent to
\begin{equation}\label{feqn4ThmA}
\Phi(w)=w,
\end{equation}
namely, $w$ is a fixed point of $\Phi$.

By Lemma~(\ref{LThmA}) bellow, we know that for any $\ve\in
(0,\rho)$ there exists $\delta=\delta(\ve)$ such that for any
$\delta$-pseudo orbit $\{x_k\}_{k\in \mathbb{Z}}$, the operator
$\Phi: \mathfrak{B}_1(\ve)\to \mathfrak{B}_1(\ve)$ defined as above
is a contracting map, and therefore has a fixed
point in $\mathfrak{B}_1(\ve)$.  Hence, (\ref{eqn4ThmA}) has a
unique solution.
\end{proof}

\begin{lemma}\label{LThmA}
For any $\ve\in
(0,\rho)$ there exists $\delta=\delta(\ve)>0$ such
that for any $\delta$-pseudo orbit
$\{x_k\}_{k\in \mathbb{Z}}$,
$\Phi(\mathfrak{B}_1(\ve))\subset \mathfrak{B}_1(\ve)$ and for any
$w,w'\in \mathfrak{B}_1(\ve)$,
$$
\|\Phi(w)-\Phi(w')\|_1 \leq
\frac{1}{2}\|w-w'\|_1.
$$
\end{lemma}

\begin{proof}
Recall that $\lambda\in (0,1)$ is given in the definition of partially hyperbolic
diffeomorphism.
For any $\widetilde{\lambda}\in (\lambda,1)$ and $\varepsilon\in (0,\rho)$,
we take $\delta>0$ and $C(\delta)>0$ such that $d(f(y),x)<\delta$ implies
\begin{equation}\label{AsThmA}
\big\|\Pi_x^s\circ d_0(\exp_x^{-1}\circ f\circ \exp_y)|_{E_y^s}\big\|
\leq \widetilde{\lambda},
\end{equation}
\begin{equation}\label{AuThmA}
\big\|\big[\Pi_x^u\circ d_0(\exp_x^{-1}\circ f\circ \exp_y)|_{E_y^u}\big]^{-1}\big\|
\leq \widetilde{\lambda},
\end{equation}
\begin{equation}\label{AuscThmA}
\sum_{i,j=s,c,u,\;i\neq j}\big\|\Pi_x^i\circ d_0(\exp_x^{-1}\circ f\circ
\exp_y)|_{E_y^j}\big\|
\leq \frac{C(\delta)}{2}
\end{equation}
and for any $v',v''\in H_x(\varepsilon)$ and any $t\in[0,1]$,
\begin{equation}\label{dbetaThmA}
\|d_{v''+t(v'-v'')}(\exp_x^{-1}\circ f\circ
\exp_y)-d_0(\exp_x^{-1}\circ
f\circ \exp_y)\|\leq\frac{C(\delta)}{2}.
\end{equation}
We can take $C(\delta)>0$ in a way such that $C(\delta)\to 0$ as $\delta\to 0$.

Note that if $\delta$ satisfies (\ref{AsThmA})--(\ref{dbetaThmA}), then
Sublemma~\ref{SL1ThmA} and \ref{SL2ThmA} below can be applied.

Further, we assume $\delta$ and $C(\delta)$ are small enough such that
\begin{equation}\label{f1LThmA}
\frac{L}{1-\tilde{\lambda}}\delta<\frac{1}{2}\varepsilon,
\qquad
\frac{L}{1-\tilde{\lambda}}C(\delta)<\frac{1}{2}.
\end{equation}

Take  $w=u+v\in \mathfrak{B}_1(\ve)$ with $u\in \mathfrak{X}^c$
and $v\in \mathfrak{X}^{us}$.
Note that for any $k\in \mathbb{Z}$,
$\|(\eta(0))_k\|=\|(\beta(0))_k\|=\|\exp^{-1}_{x_k}f(x_{k-1})\|\le \delta$.
and hence $\|\eta(0)\|\le \delta$.
So by Sublemma~\ref{SL1ThmA} and \ref{SL2ThmA} below, and \eqref{f1LThmA},
we can get
\begin{equation*}
\begin{split}
&\|\Phi(w)\|_1
\leq \|P^{-1}\|_1\cdot\|\eta(v)\|_1
\leq \frac{1}{1-\tilde{\lambda}}\cdot L  \|\eta(v)\| \\
\leq &\frac{L}{1-\tilde{\lambda}}
   (\|\eta(v)-\eta(0)\|+\|\eta(0)\|)
\leq \frac{L}{1-\tilde{\lambda}}
   (C(\delta)\|v\|_1+\delta)
< \frac{1}{2}\|w\|_1
    +\frac{1}{2}\varepsilon
\leq \varepsilon,
\end{split}
\end{equation*}
which implies that
$\Phi(\mathfrak{B}_1(\ve))\subset\mathfrak{B}_1(\ve)$.

Similarly, for two elements $w=u+v,\;w'=u'+v'\in
\mathfrak{B}_1(\ve)$ with $u,u'\in \mathfrak{X}^c$ and $v,v'\in
\mathfrak{X}^{us}$, we have
\begin{equation*}
\begin{split}
&\big\|\Phi(w)-\Phi(w')\big\|_1
\leq\frac{1}{1-\tilde{\lambda}}  \big(\|\eta(v)-\eta(v')\|_1\big) \\
\leq &\frac{L }{1-\tilde{\lambda}}\big(\|\eta(v)-\eta(v')\|\big)
\leq \frac{L }{1-\tilde{\lambda}}\big(C(\delta)\|w-w'\|_1\big)
\leq\frac{1}{2}\|w-w'\|_1.
\end{split}
\end{equation*}
This proves that $\Phi: \mathfrak{B}_1(\ve)\to
\mathfrak{B}_1(\ve)$ is a contraction.
\end{proof}

\begin{sublemma}\label{SL1ThmA}
For $\delta>0$ satisfying (\ref{AsThmA})--(\ref{dbetaThmA})
and any $v, v'\in \mathfrak{B}^{us}(\ve)$,
\begin{equation*}\label{equation3 of lemma01}
\|\eta(v')-\eta(v)\|\leq C(\delta)(\|v'-v\|),
\end{equation*}
where $C(\delta)$ is chosen in the beginning of the proof of Lemma~\ref{LThmA}.
\end{sublemma}

\begin{proof}
Denote $\eta_{k}(v_k)=(\eta(v))_{k+1}$ for $v=\{v_k\}_{k\in \mathbb{Z}}$
in a neighborhood of $0\in \mathfrak{X}^{us}$.
By the definition of $\eta$, we can write
$$
\eta_k=\eta_k^{(1)}+\eta_k^{(2)},
$$
where
$$
\eta_k^{(1)}(v_k)=\exp_{x_{k+1}}^{-1}\circ f\circ
\exp_{x_k}(v_k)-d_0(\exp_{x_{k+1}}^{-1}\circ f\circ
\exp_{x_k})v_k
$$
and
$$
\eta_k^{(2)}(v_k)=\sum_{i=s,c,u,\;j=s,u,\;i\neq j}\Pi_{x_{k+1}}^i\circ
d_0(\exp_{x_{k+1}}^{-1}\circ f\circ
\exp_{x_k})\circ \Pi_{x_k}^jv_k.
$$
Note that for $v',v''\in H_k(\varepsilon)$, we have
\begin{eqnarray}
&&\big\|\eta_k^{(1)}(v')-\eta_k^{(1)}(v'')\big\|\notag\\
&=&\Big\|\int_0^1 \big[d_{v''+t(v'-v'')}(\exp_{x_{k+1}}^{-1}\circ f\circ
\exp_{x_k})-d_0(\exp_{x_{k+1}}^{-1}\circ
f\circ \exp_{x_k})\big](v'-v'')dt\Big\|\notag\\
&\leq& \sup_{t\in[0,1]}\big\|d_{v''+t(v'-v'')}(\exp_{x_{k+1}}^{-1}\circ f\circ
\exp_{x_k})-d_0(\exp_{x_{k+1}}^{-1}\circ f\circ \exp_{x_k})\big\|\cdot
\big\|v'-v''\big\|.\notag
\end{eqnarray}
Therefore, from (\ref{dbetaThmA}) we have
\begin{equation}\label{omega1}
\|\eta_k^{(1)}(v')-\eta_k^{(1)}(v'')\|\leq \frac{C(\delta)}{2}
\|v'-v''\|.
\end{equation}

By  (\ref{AuscThmA}), we have for $v',v''\in H_k(\varepsilon)$
\begin{equation}\label{omega2}
\|\eta_k^{(2)}(v')-\eta_k^{(2)}(v'')\|\leq\frac{C(\delta)}{2}\|v'-v''\|.
\end{equation}

Combining (\ref{omega1}) and (\ref{omega2}), for
$v',v''\in H_k(\varepsilon)$ we have
\begin{equation}\label{omega}
\|\eta_k(v')-\eta_k(v'')\|\leq C(\delta) \|v'-v''\|.
\end{equation}
Hence, we can get the result we need immediately.
\end{proof}

\begin{sublemma}\label{SL2ThmA}
For any $\delta>0$ satisfying (\ref{AsThmA})--(\ref{dbetaThmA})
and any $\delta$-pseudo orbit $\{x_k\}_{k\in \mathbb{Z}}$, the operator $P$
defined as (\ref{PThmA}) is invertible and
$$
\|P^{-1}\|_1\leq \frac{1}{1-\tilde{\lambda}}.
$$
\end{sublemma}

\begin{proof}
By the definition of $P$, we have
$P|_{\mathfrak{X}^i}=id_{\mathfrak{X}^i}-A^i,\;i=s,u$,
and $P|_{\mathfrak{X}^c}=\id_{\mathfrak{X}^c}$. So
$P(\mathfrak{X}^i)=\mathfrak{X}^i,\;i=u,s,c$.

By (\ref{AuThmA}) and (\ref{AsThmA}), $\|A^s\|, \|(A^u)^{-1}\|\le \tilde{\lambda}<1$.
Hence, both $P|_{\mathfrak{X}^s}$ and $P|_{\mathfrak{X}^u}$ are
invertible and
\begin{equation*}
\begin{split}
(P|_{\mathfrak{X}^s})^{-1}
&=(id_{\mathfrak{X}^s}-A^s)^{-1}
=\sum_{k=0}^\infty A^s_k,      \\
(P|_{\mathfrak{X}^u})^{-1}
&=(id_{\mathfrak{X}^u}-A^u)^{-1}=-\sum_{k=1}^\infty (A^u_k)^{-1}.
\end{split}
\end{equation*}
It follows that
\begin{eqnarray*}
\|(P|_{\mathfrak{X}^{us}})^{-1}\| \le
\max\left\{\|(P|_{\mathfrak{X}^s})^{-1}\|,
               \|(P|_{\mathfrak{X}^u})^{-1}\|\right\}
\le \frac{1}{1-\tilde{\lambda}}.
\end{eqnarray*}
It is obvious that
\begin{equation*}
\|(P|_{\mathfrak{X}^c})^{-1}\|=1.
\end{equation*}
So we obtain that
$$
\|P^{-1}\|_1 \le \max\left\{\|(P|_{\mathfrak{X}^{us}})^{-1}\|,
               \|(P|_{\mathfrak{X}^c})^{-1}\|\right\}
\le \frac{1}{1-\tilde{\lambda}}.
$$
This is what we need.
\end{proof}

\section{Quasi-shadowing for the system with $C^1$ center foliation}

\subsection{The general case}

Recall that $\mathfrak{X}^{us}$ and $\mathfrak{B}^{us}(\rho)$ are defined
in the beginning of the previous section.

\begin{proof}[Proof of Theorem B]
The proof is similar to that of Theorem A.

Given a $\delta$-pseudo orbit $\{x_k\}_{k\in \mathbb{Z}}$ of
$f$. To find a sequence of points $\{y_k\}_{k\in \mathbb{Z}}$ satisfying (\ref{eqn1thmB}), (\ref{eqn2thmB}) and (\ref{eqn3thmA}), we shall
try to solve the equations
\begin{equation}\label{feqn2ThmB}
y_k=\tau_{x_k}^{(2)}(f(y_{k-1}))
\end{equation}
for unknown $\{y_k\}_{k\in \mathbb{Z}}$.
Put $v_k=\exp_{x_k}^{-1}y_k$.
Then the equations (\ref{feqn2ThmB}) are equivalent to
\begin{equation}\label{eqn4ThmB}
v_{k}=\exp_{x_{k}}^{-1}\circ \tau_{x_{k}}^{(2)}\circ f\circ \exp_{x_{k-1}}v_{k-1},
\qquad  k\in \mathbb{Z}.
\end{equation}

Define an operator $\beta: \mathfrak{B}^{us}(\rho)\to \mathfrak{X}^{us}$
and a linear operator $A: \mathfrak{B}^{us}(\rho)\to \mathfrak{X}^{us}$ by
\begin{equation}\label{betaThmB}
(\beta(v))_k
=\exp_{x_k}^{-1}\circ \tau_{x_k}^{(2)}\circ f\circ \exp_{x_{k-1}}v_{k-1},
\quad k\in \mathbb{Z},
\end{equation}
and
\begin{equation}\label{AThmB}
(Av)_k
=(A_{k-1}^u+A_{k-1}^s)v_{k-1},
\end{equation}
where
\begin{equation*}
\begin{split}
A_{k-1}^s=\Pi_{x_k}^s\circ
d_0(\exp_{x_k}^{-1}\circ\tau_{x_k}^{(2)}\circ f\circ \exp_{x_{k-1}})\circ \Pi_{x_{k-1}}^s, \\
A_{k-1}^u=\Pi_{x_k}^u\circ
d_0(\exp_{x_k}^{-1}\circ\tau_{x_k}^{(2)}\circ f\circ \exp_{x_{k-1}})\circ \Pi_{x_{k-1}}^u.
\end{split}
\end{equation*}

Let $\eta=\beta-A$.  By (\ref{betaThmB}) and (\ref{AThmB}),
(\ref{eqn4ThmB}) is equivalent to
$$
v=Av+\eta(v),
$$
further, is equivalent to
$$
v-Av=\eta(v).
$$

Define a linear operator $P$ from a neighborhood of $0\in
\mathfrak{X}^{us}$ to $\mathfrak{X}^{us}$ by
\begin{equation}\label{PThmB}
Pv=(id_{\mathfrak{X}^{us}}-A)v,
\end{equation}
and then define an operator $\Phi$ from a neighborhood of $0\in
\mathfrak{X}^{us}$ to $\mathfrak{X}^{us}$ by
\begin{equation}\label{DefPhi}
\Phi(v)=P^{-1}\eta(v)
\end{equation}
for $v$ in a neighborhood of $0\in
\mathfrak{X}^{us}$.

Hence, the equations (\ref{eqn4ThmB}) are equivalent to
\begin{equation}\label{feqn4ThmA}
\Phi(v)=v,
\end{equation}
namely, $v$ is a fixed point of $\Phi$.

The remaining work is to show that for any $\ve\in
(0,\rho)$ there exists $\delta=\delta(\ve)$ such that
for a $\delta$-pseudo orbit $\{x_k\}_{k\in \mathbb{Z}}$ of
$f$, $\Phi: \mathfrak{B}^{us}(\ve)\to
\mathfrak{B}^{us}(\ve)$ is a contracting map, and therefore has a fixed
point in $\mathfrak{B}^{us}(\ve)$. Hence, (\ref{eqn4ThmB}) has a
unique solution. To this end we only need to modify the proof
of Lemma~\ref{LThmA} to a easer version
since in this case we do not need to consider the center direction.
We leave the details to the reader.
\end{proof}

\subsection{$\mathcal{W}^c_f$ is of one dimensional}

\begin{proof}[Proof of Theorem B$'$]
The proof is also similar to that of Theorem A.

Given a $\delta$-pseudo orbit $\{x_k\}_{k\in \mathbb{Z}}$ of
$f$. To find a sequence of points $\{y_k\}_{k\in \mathbb{Z}}$ and
and a sequence of real numbers $\{\tilde{\tau}_k\}_{k\in \mathbb{Z}}\in \mathfrak{C}(\rho)$, where
$\mathfrak{C}(\rho)=\{\tilde{\tau}=\{\tilde{\tau}_k\}_{k\in \mathbb{Z}}:\tilde{\tau}_k\in \mathbb{R}, |\tilde{\tau}_k|\le \rho, k\in \mathbb{Z}
\}$, satisfying (\ref{eqn1thmB'}), (\ref{eqn2thmB'}) and (\ref{eqn3thmA}), we shall
try to solve the equations
\begin{equation}\label{feqn2ThmB'}
y_k=\tau_{x_k}^{(3)}(f(y_{k-1}))
\end{equation}
for unknown $\{y_k\}_{k\in \mathbb{Z}}$ and $\{\tilde{\tau}_k\}_{k\in \mathbb{Z}}$.
Putting $v_k=\exp_{x_k}^{-1}y_k$, then the equations (\ref{feqn2ThmB'}) are equivalent to
$$
v_{k+1}=\exp_{x_{k+1}}^{-1}\circ \tau_{x_{k+1}}^{(3)}\circ f\circ \exp_{x_k}v_k,\;\;\;k\in \mathbb{Z},
$$
i.e.,
\begin{equation}\label{eqn4ThmB'}
v_{k+1}=\exp_{x_{k+1}}^{-1}\circ \varphi^{\tilde{\tau}_{k+1}}\circ f\circ \exp_{x_k}v_k,\;\;\;k\in \mathbb{Z}.
\end{equation}

Define an operator $\beta: \mathfrak{B}^{us}(\rho)\times \mathfrak{C}(\rho)\to
\mathfrak{X}$ and a linear operator $A: \mathfrak{B}^{us}(\rho)\to \mathfrak{X}^{us}$ by
\begin{eqnarray}\label{betaThmB'}
\beta(v,\tilde{\tau})_{k+1}=\exp_{x_{k+1}}^{-1}\circ
\varphi^{\tilde{\tau}_{k+1}} \circ f
  \circ \exp_{x_k}\bigl(v_k)
\end{eqnarray}
and
\begin{equation}\label{AThmB'}
(Av)_k
=(A_{k-1}^u+A_{k-1}^s)v_{k-1},
\end{equation}
where
\begin{equation*}
\begin{split}
A_{k-1}^s=\Pi_{x_k}^s\circ
(d_{(0,0)}\beta(v,\tilde{\tau}))_k\circ \Pi_{x_{k-1}}^s,  \\
A_{k-1}^u=\Pi_{x_k}^u\circ(d_{(0,0)}\beta(v,\tilde{\tau}))_k\circ \Pi_{x_{k-1}}^u.
\end{split}
\end{equation*}
Let $\eta=\beta-A$.
Let $u$ be a vector field consisting of unit vectors tangent to $\mathcal{W}^c_f$.
Then by (\ref{betaThmB'}) and (\ref{AThmB'}), (\ref{eqn4ThmB'}) is equivalent to
$$
v=\tilde{\tau}u+Av+\eta(v)
$$
for some $\tilde{\tau}\in \mathfrak{C}(\rho)$.
Further, the equations are equivalent to
$$
-\tilde{\tau}u+v-Av=\eta(v).
$$

Define a linear operator $P$ from a neighborhood of $0\in \mathfrak{X}(\rho)$
to $\mathfrak{X}$ by
\begin{equation}\label{PThmB'}
Pv=\tilde{\tau}u+(id_{\mathfrak{X}^{us}}-A)v,
\end{equation}
and then define an operator $\Phi$ from a neighborhood of $0\in
\mathfrak{X}$ to $\mathfrak{X}$ by
$$
\Phi(w)=P^{-1}\eta(v),
$$
where $w=\tilde{\tau}\cdot u+v\in \mathfrak{X}$ with
$v\in \mathfrak{B}^{us}(\rho)$ and $\tilde{\tau}\in \mathfrak{C}(\rho)$.

Hence, the equations (\ref{eqn4ThmB'}) are equivalent to
\begin{equation}\label{feqn4ThmB'}
\Phi(\tilde{\tau}\cdot u+v)=\tilde{\tau}\cdot u+v,
\end{equation}
namely, $\tilde{\tau}\cdot u+v$ is a fixed point of $\Phi$.

The remaining work is to show that for any $\ve\in
(0,\rho)$ there exists $\delta=\delta(\ve)$ such that for $\delta$-pseudo orbit $\{x_k\}_{k\in \mathbb{Z}}$ of
$f$, $\Phi: \mathfrak{B}_1(\varepsilon)\to
\mathfrak{B}_1(\ve)$ is a contracting map, and therefore has a fixed
point in $\mathfrak{B}_1(\varepsilon)$. Hence, (\ref{eqn4ThmB'}) has a
unique solution. We leave the details to the reader.
\end{proof}



\section{Quasi-stability}

It is well known that for any homeomorphism $f$ on a compact metric space,
shadowing property together with expansiveness implies topological stability
(see [16] for example).
In the case of partially hyperbolic diffeomorphism,
we can get topological quasi-stability from quasi-shadowing property.

\begin{proof}[Proof of Theorem C]

For simplicity of the notation, we only prove this theorem
under the condition that $f$ is a partially hyperbolic diffeomorphism
with $C^1$ center foliation.

Choose $\varepsilon_0\in (0,\rho)$ small enough such that
any continuous map $h$ with $d(h,\id_M)<\varepsilon_0$ must be
surjective (see e.g Lemma 3 of \cite{Walters}
for existence of such $\ve_0$).

Let $\varepsilon\in (0, \varepsilon_0)$.
From Theorem B, there exists $\delta>0$ such that for any $\delta$-pseudo
orbit $\{x_k\}_{k\in \mathbb{Z}}$ of $f$, there exists a unique pseudo orbit
$\{y_k\}_{k\in \mathbb{Z}}$ $\varepsilon$-quasi-shadowing it that satisfies
(\ref{eqn3thmA}) and $y_{k+1}\in W^c(y_k)$ for all $k\in \mathbb Z$.
Let $g$ be any homeomorphism with $d(f,g)<\delta$.
It is obvious that for any $x\in M$, its orbit
$\textmd{orb}_g(x)=\{x_k=g^k(x)\}_{k\in \mathbb{Z}}$
is a $\delta$-pseudo orbit of $f$, hence, there exists a unique
corresponding pseudo orbit $\{y_k\}_{k\in \mathbb{Z}}$
$\varepsilon$-quasi-shadowing it. Let $h(x)=y_0$.

Now we consider continuity of $h$.
Notice that the sequence $\{y_k\}_{k\in \mathbb{Z}}$,
which is $\varepsilon$-quasi-shadowing the orbit of $x$,
is defined by the sequence $v=\{\exp^{-1}_{g^k(x)}y_k\}_{k\in \mathbb{Z}}$,
and $v$ is the fixed point of the operator $\Phi_{\textmd{orb}_g(x)}:
\mathfrak{B}_{\textmd{orb}_g(x)}^{us}(\ve)\to \mathfrak{B}_{\textmd{orb}_g(x)}^{us}(\ve)$
in the proof of Theorem B
(here we use the notions $\Phi_{\textmd{orb}_g(x)}$ and
$\mathfrak{B}_{\textmd{orb}_g(x)}^{us}(\ve)$
since they all depend on $\textmd{orb}_g(x)$).
Given $x'$ near $x$, denote by $v'=\{v'_k\in E^{us}_{g^k(x')}\}$
the unique fixed point of the operator $\Phi_{\textmd{orb}_g(x')}:
\mathfrak{B}_{\textmd{orb}_g(x')}^{us}(\ve)\to \mathfrak{B}_{\textmd{orb}_g(x')}^{us}(\ve)$. By the definition of $h$, $h(x')=\exp_{x'}(v'_0)$.
By continuity of the distribution $E^{us}$, continuity of the differential of $f$ and the construction of the operator $\Phi$, we can see that as $x'$ approaches $x$, $v'_0$ approaches $v_0$ in the tangent bundle $TM$. Therefore, $h(x')$ arbitrarily approaches $h(x)$ as $x'$ sufficiently close to $x$. This means that the map $h$ is continuous.
\end{proof}

\section{Center Nonwandering Sets}

It is well known that for a uniformly hyperbolic system $f: M\to M$,
if $x$ is close to $f^nx$ for some $x\in M$ and $n>0$,
then there is a periodic point $y\in M$ of period $n$ close to $x$.
The result is the main part of Anosov closing lemma
(see e.g. \cite{Bowen, Shub}),
and sometimes is directly called Anosov closing lemma
(see e.g. \cite{KH}).

The next lemma is an analogue of the result for
partially hyperbolic diffeomorphisms.

\begin{lemma}\label{L1ThmD}
Suppose $f: M\to M$ is a partially hyperbolic diffeomorphism
with $C^1$ center foliation $\mathcal{W}^c_f$.
For any $\varepsilon>0$, there exists $\delta\in (0,\varepsilon)$ such that
for any $x\in M$ and $n\in {\mathbb N}$ with $d(x, f^nx)<\delta$,
there exists a periodic center leaf $W^c(p)$ of period $n$ satisfying
$d(p, x)\le \varepsilon$.

Moreover, if $W^c(x)$ is compact and
$d_H\bigl(W^c(x),f^n(W^c(x))\bigr)<\delta$,
then there exists a periodic center leaf $W^c(p)$ of period $n$
such that $d(p, x')\le \varepsilon$ for some $x'\in W^c(x)$.
\end{lemma}

\begin{proof}
By Theorem B, there is $\delta\in (0,\varepsilon)$
such that any $\delta$-pseudo orbit can be $\varepsilon$-quasi-shadowed.
Since $d(x, f^nx)<\delta$, we can repeat the orbit segment
$\{x, fx, \cdots, f^{n-1}x\}$ to get a $\delta$-pseudo orbit
$\{x_k\}_{k\in \mathbb{Z}}$, where
$x_k=f^ix$ if $k\equiv i (\mod n)$.
By Theorem B, there is a unique sequence $\{y_k\}_{k\in \mathbb{Z}}$ which
$\varepsilon$-quasi-shadows $\{x_k\}_{k\in \mathbb{Z}}$ and satisfies (\ref{eqn3thmA}) and (\ref{eqn2thmB}).
Note that  $x_{n+i}=x_i$ for all $i\in \mathbb{Z}$,
$\{y_{n+k}\}_{k\in \mathbb{Z}}$ also $\varepsilon$-quasi-shadows
$\{x_{k}\}_{k\in \mathbb{Z}}$.

By the uniqueness of the quasi-shadowing sequence $\{y_k\}_{k\in \mathbb{Z}}$ satisfying (\ref{eqn3thmA}) and (\ref{eqn2thmB}), $y_n=y_0$.
Note that if $\{y_k\}_{k\in \mathbb{Z}}$ $\varepsilon$-quasi-shadows
$\{x_k\}_{k\in \mathbb{Z}}$, then $f(y_{k-1})\in W^c(y_k)$ and therefore
$fW^c(y_{k-1})=W^c(y_k)$.
So $y_n=y_0$ implies $f^nW^c(y_0)=W^c(y_n)=W^c(y_0)$,
that is, $W^c(y_0)$ is a periodic center leaf.
So the first part of the lemma follows with $p=y_0$.

To prove the second part, we use the inequality
$d_H\bigl(W^c(x),W^c(f^nx)\bigr)<{\delta}$.
Take $x_0=x$ and then take $x_1\in W^c(x)$ such that $d(x_1, f^n(x_0))<\delta$.
Inductively, for any $i\ge 1$, if $x_{i-1}\in W^c(x)$ is taken,
then we can choose $x_{i}\in W^c(x)$ such that $d(f^n(x_{i-1}), x_i)< \delta$.
Since $W^c(x)$ is compact, there exist $i<j$ such that $d(x_i, x_j)<\delta-d(f^n(x_{j-1}), x_j)$.
Repeating the pseudo orbit segment
$$x_i, f(x_i), \cdots, f^{n-1}(x_i), x_{i+1}, f(x_{i+1}),
\cdots, f^{n-1}(x_{i+1}),\cdots, x_{j-1}, f(x_{j-1}), \cdots, f^{n-1}(x_{j-1})$$
we get a $\delta$-pseudo orbit $\{x^{(k)}\}_{k\in \mathbb{Z}}$
satisfying that $x^{(0)}=x_i$ and $x^{(k)}=x^{(\ell)}$
if $k\equiv \ell$ $(\mod n(j-i))$.

By the same arguments as above, there is a periodic
center leaf $W^c(y)$ with $d(x^{(0)},y)\le \varepsilon$.
Hence, we get the result if we take $p=y$ and $x'=x^{(0)}$.
\end{proof}

\begin{proof}[Proof of Theorem D]
The first result of the theorem follows from
the first result of Lemma~\ref{L1ThmD} immediately.
This is because for any $x\in \Omega(f)$ and $\varepsilon>0$,
we can find $y\in M$ and $n>0$ such that
$d(x,y)\le \varepsilon$ and $d(y,f^ny)\le \delta$ (here $\delta>0$ is as that in Lemma ~\ref{L1ThmD}),
and therefore there exists $p\in P^c(f)$ with $d(y, p)<\varepsilon$.
Hence, $d(x,p)\le 2\varepsilon$.
It means that $x\in \overline{P^c(f)}$.
We get $\Omega(f)\subset \overline{P^c(f)}$.

Now we consider the second part of the theorem.
For any $\varepsilon>0$, choose $\delta>0$ as in Lemma ~\ref{L1ThmD}. Note that if $f$ has uniformly compact center foliation, then
by (\ref{fOmegac}), for any $x\in \Omega^c(f)$ there is
$y\in M$ such that $d_H\bigl(W^c(x),W^c(y)\bigr)<\varepsilon$
and $d_H\bigl(W^c(y),f^n(W^c(y))\bigr)<{\delta}$.
By the same arguments we get that
 there exists $p\in P^c(f)$ with $d(y', p)<\varepsilon$
for some $y'\in W^c(y)$.
Hence $d(x', p)<2\varepsilon$ for some $x'\in W^c(x)$.
And we get $\Omega^c(f)\subset\overline{P^c(f)}$.

By Lemma~\ref{L2ThmD} below, we know that
$\Omega^c(f)\subset W^c(\Omega(f))$.
Since it is obvious that $\overline{P^c(f)}\subset\Omega^c(f)$ and
$W^c(\Omega(f))\subset\Omega^c(f)$, we get the equality (\ref{fThmD}).
\end{proof}

\begin{lemma}\label{L2ThmD}
Suppose $f: M\to M$ is a partially hyperbolic diffeomorphism
with compact $C^1$ center foliation $\mathcal{W}^c_f$.
Then for any $x\in \Omega^c(f)$, there exists $x'\in \Omega(f)\cap W^c(x)$.
\end{lemma}

\begin{proof}
Suppose $\Omega(f)\cap W^c(x)=\emptyset$.  Then any point $y\in W^c(x)$ is
a wandering point.  Hence, there is a neighborhood $U_y$ of $y$ such that
$f^n(U_y)\cap U_y=\emptyset$ for any $n>0$.
Clearly $\{U_y: y\in W^c(x)\}$ form a open cover of $W^c(x)$.
Let $\{U_1, \dots, U_k\}$ be a subcover of $W^c(x)$,
and let $U=\cup_{i=1}^k U_i$.  Then $U$ is a neighborhood of $W^c(x)$.
By Lemma~\ref{L1ThmD}, $U$ contains a periodic leaf $W^c(z)$.

Suppose $W^c(z)$ has period $\ell$.  Then $f^{j\ell}(z)\in U$ for any $j>1$.
Since $U=\cup_{i=1}^k U_i$, there are $j_1<j_2$ such that
$f^{j_1 \ell}(z), f^{j_2 \ell}(z)\in U_i$ for some $U_i$.
That is, $f^{(j_2-j_1)\ell}U_i\cap U_i \not=\emptyset$,
contradicting the fact that $f^n(U_y)\cap U_y=\emptyset$ for any $y\in W^c(x)$ and any $n>0$.
\end{proof}

\section{A spectral decomposition theorem}

In this section, we assume that $f: M\to M$ is a partially hyperbolic
diffeomorphism with uniformly compact $C^1$ center foliation.

Denote by $W^u_\varepsilon(x)$ and $W^s_\varepsilon(x)$ the local unstable
and stable manifolds of size $\varepsilon$ at $x$ respectively.

We recall that $f$ is dynamically coherent since the center foliation
$W^c$ of $f$ is $C^1$ (\cite{Pugh}).

The next lemma gives the local product structure.

\begin{lemma}\label{L0ThmE}
There are $\varepsilon, \delta>0$ such that for any $x, y\in M$ with
$d(x, y)<\delta$, for any $x_1\in W^c(x)$, there is $y_1\in W^c(y)$
such that $W^s_\ve(x_1)\cap W^u_\ve(y_1)$ contains exact one point.
\end{lemma}

\begin{proof}
Since $W^s$ and $W^{cu}$ are uniformly transversal,
and $W^u$ subfoliate $W^{uc}$,
it is obvious that there are $\varepsilon, \delta'>0$ such that
if $w, z\in M$ with $d(w, z)<\delta'$,
then $W^s_\ve(w)\cap W^u_\ve(z_1)$ contains exact one point
for some $z_1\in W^u(z)$.

Since $f$ has uniformly compact $C^1$ center foliation,
we can take $\delta\in (0,\delta')$ such that if $d(x,y)< \delta$, then
$d_H\bigl(W^c(x),W^c(y)\bigr)<\delta'$.
Then for any $x_1\in W^c(x)$, we can find $y_0\in W^c(y)$
such that $d(x_1, y_0)<\delta'$.  Thus the result follows.
\end{proof}

For uniformly hyperbolic systems, the cloud lemma gives that
for any periodic points $p$ and $q$, any point $x\in W^u(p)\cap W^s(q)$
is contained in the nonwandering set of the map (see e.g. \cite{Shub}).
The next lemma can be regarded as a local version
of the cloud lemma for partially hyperbolic diffeomorphisms
with uniformly compact $C^1$ center foliation.

\begin{lemma}\label{CloudLemma}
Let $\delta$ and $\varepsilon$ be positive
numbers as in Lemma \ref{L0ThmE} and $p, q\in P^c(f)$ with $d(p, q)< \delta$.
If $x\in W^s_\ve(p_1)\cap W^u_\ve(q_1)$ for some $p_1\in W^c(p)$ and
$q_1\in W^c(q)$, then $x\in\Omega^c(f)$.
\end{lemma}

\begin{figure}[htb]
\begin{center}
\psfrag{P}{$p_1$} \psfrag{Q}{$q$} \psfrag{X}{$x$}
\psfrag{A}{$W^c(p)$} \psfrag{B}{$W^c(x)$}\psfrag{K}{$W^c(q)$}
\psfrag{C}{$E^c$}\psfrag{U}{$E^u$}\psfrag{S}{$E^s$}
\psfrag{Y}{$y$} \psfrag{H}{$z$} \psfrag{1}{$q_1$} \psfrag{2}{$q_2$}\psfrag{3}{$q_3$}
\psfrag{G}{$f^{i_0}x$}\psfrag{D}{$f^{j_0}z$}
\psfrag{J}{$p_2$}\psfrag{L}{$W^s(p_1)$}\psfrag{M}{$W^u(q_1)$}
\psfrag{W}{$z_0$}\psfrag{E}{$f^{j_0}y$}\psfrag{F}{$x'$}
\psfrag{N}{$p$}
\includegraphics{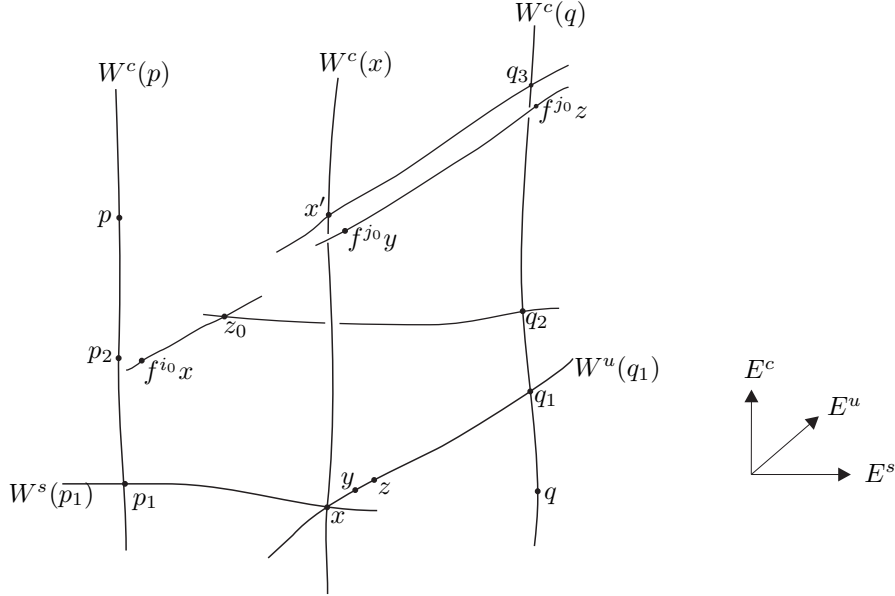}
 \caption{intersection of stable and unstable manifolds} \label{figure}
\end{center}
\end{figure}

\begin{proof}
By the definition of $\Omega^c(f)$ and uniform compactness
of the center foliation, it is sufficient to prove that for any
$\alpha>0$, there are a point $y$ and a number $n\in \mathbb{N}$
such that
\begin{equation}\label{eq2ThmE}
d(x,y)<\alpha  \qquad \text{ and }\qquad d(f^n(y),W^c(y))<\alpha.
\end{equation}
By uniform compactness of the center foliation, there exists
$\beta\in (0,\alpha)$ such that
\begin{equation}\label{eq0ThmE}
d(x,y)<\beta \
\quad \Longrightarrow \quad
d_H(W^c(x), W^c(y)) <\frac{\alpha}{2} \qquad \forall x,y\in M.
\end{equation}

Since $x\in W^s(p_1)$ and $p_1\in P^c(f)$,
$\{f^i(x)\}_{i\geq 0}$ has an accumulation point $p_2\in W^c(p)$
(See Figure~\ref{figure}).
Note that $d(p_2, W^c(q))< \delta'$ by the choice of $\delta$
in the proof of Lemma~\ref{L0ThmE}.
Hence, we can find $i_0>0$ such that $d(f^{i_0}(x), W^c(q))< \delta'$.
By Lemma~\ref{L0ThmE}, there are $q_2\in W^c(q)$ and $z_0\in M$
such that $z_0\in W^u_\ve(f^{i_0}(x))\cap W^s_\ve(q_2)$.
Set $z=f^{-i_0}(z_0)$.  We can choose  $i_0$ large enough such that
$$d(z,x)<\frac{\beta}{2},
$$
where $\beta$ is given in (\ref{eq0ThmE}).
Note that $\{f^n(z_0)\}_{n\geq 0}$ has an accumulation point $q_3\in W^c(q)$
since $z\in W^{cs}(q_2)=W^{cs}(q)$ and $W^c(q)$ is compact.
We may assume $f^{n_j}(z_0)\to q_3$ for some $n_j\to +\infty$.
This implies that
\begin{equation}\label{eq1ThmE}
\lim\limits_{j\to \infty} f^{n_j-i_0}(z)=q_3.
\end{equation}

Recall that $x\in W^{u}_\ve(q_1)\subset W^{cu}(q)$.
There is a point $x'\in W^c(x)\cap W^u(q_3)$ since $W^c$ and $W^u$ subfoliate
$W^{cu}$.

Note that $z\in W^u_\ve(x)$ and $x'\in W^u(q_3)$.
By continuity of the unstable foliation, (\ref{eq1ThmE})
implies that there are a point $y\in W^u_{\frac{\beta}{2}}(z)$ and
an integer $j_0\in \mathbb{N}$ such that
$d(f^{j_0}(y),x')<\frac{\alpha}{2}$,
and therefore
$$
d(f^{j_0}(y),W^c(x))<\frac{\alpha}{2}.
$$
Also we have that $d(x,y)\le d(x,z)+d(z,y)
< \frac{\beta}{2}+\frac{\beta}{2}=\beta$.
So (\ref{eq0ThmE}) can be applied, and we have
$$
d(f^{j_0}(y),W^c(y))\le d(f^{j_0}(y),W^c(x))+ d(W^c(x), W^c(y))
<\frac{\alpha}{2}+\frac{\alpha}{2}
=\alpha.
$$
Now we get (\ref{eq2ThmE}) with $n=j_0$, and complete the proof.
\end{proof}

\begin{proof}[Proof of Theorem E]
For any $p\in P^c(f)$, we set
$$
X_p=\overline{W^{cu}(p)\cap \Omega^c(f)}.
$$
By Lemma~\ref{L1ThmE} below, $X_p$ is both open and closed in
$\Omega^c(f)$.
By Lemma~\ref{L2ThmE} below, any two elements in $\{X_p: p\in P^c(f)\}$
are either disjoint or equal.
Since $\Omega^c(f)$ is compact, there are finitely many points
$p_1, \cdots, p_n\in P^c(f)$ such that
$$\Omega^c(f)=X_{p_1}\cup X_{p_2}\cup \cdots \cup X_{p_n},
$$
where $X_{p_i}$ are pairwise disjoint.
Then $f(X_{p_j})=X_{f(p_j)}$ and hence equals to some $X_{p_i}$.
So $f$ permutes these $X_{p_j}$'s.
We set $\Omega^c_i$ as the union of the $X_{p_j}$'s
in the various cycles of permutation.
Then we can get
$$\Omega^c(f)=\Omega^c_1\cup\cdots\cup\Omega^c_k.
$$

Center topological transitivity  in (a) is implied
by center topological mixing in (b).
For finishing the proof, we only need to prove that
$f^N:X_p\to X_p$ is center mixing whenever $p\in P^c(f)$ and
$N\in \mathbb{N}$ satisfying $f^N(X_p)=X_p$.

Suppose $U,V$ are nonempty open subsets in $X_p$.
We choose a point $q\in P^c(f)\cap U$, and assume that $W^c(q)$ has period $n$.
Note that $n=tN$ for some $t\geq 1$.

We firstly prove that there is $i_0\in \mathbb{N}$ such that
\begin{equation}\label{eq3ThmE}
f^{in}(W^c(U))\cap V\neq \emptyset \qquad \forall i\geq i_0.
\end{equation}
In fact, since $U$ is open, there exists $\varepsilon>0$ such that
$$B_{\varepsilon,\Omega^c(f)}(W^c(q))
=\{x\in\Omega^c(f): d(x, W^c(q))<\varepsilon\}\subset W^c(U).
$$
On the other hand, since $W^{cu}(q)$ is dense in $X_q=X_p$,
we can select a point $z\in W^{cu}(q)\cap \Omega^c(f)\cap V$.
Then there is $i_0\in \mathbb{N}$ such that
$$f^{-in}(z)\in B_{\varepsilon,\Omega^c(f)}(W^c(q))\qquad \forall i\geq i_0
$$
and hence this proves (\ref{eq3ThmE}).

Similar to (\ref{eq3ThmE}), for any $j=1,\cdots,t-1$,
there is $i_j\in \mathbb{N}$ such that
\begin{equation}\label{eq4ThmE}
f^{in}(f^{jN}(W^c(U)))\cap V\neq \emptyset \qquad  \forall i\geq i_j.
\end{equation}
Set $i_*=t\cdot \max\{i_0, i_1,\cdots,i_{t-1}\}$.
Then, for any $i\geq i_*$, we can write $$iN=ln+jN,$$
where $l\geq \max\{i_0, i_1,\cdots,i_{t-1}\}$
and $1\leq j\leq t-1$. So, by (\ref{eq3ThmE}) and (\ref{eq4ThmE})
$$
f^{iN}(W^c(U))\cap V =f^{ln}(f^{jN}(W^c(U)))\cap V\neq \emptyset
\qquad  \forall i\geq i_*.
$$
That is to say, $f^N|_{X_p}$ is center mixing.
We complete the proof of Theorem E.
\end{proof}

\begin{lemma}\label{L1ThmE}
There exists $\delta>0$ such that for any $p\in P^c(f)$,
$$B_{\delta,\Omega^c(f)}(X_p):=\{x\in \Omega^c(f): \ d(x,X_p)<\delta\}=X_p,$$
where $d(x,X_p)=\min\limits_{y\in X_p}d(x,y)$.
\end{lemma}

\begin{proof}
Let $\delta>0$ be as in Lemma~\ref{L0ThmE}.
By Theorem~D, we only need to prove that $q\in X_p$
for any $q\in P^c(f)$ with $d(q, X_p)<\delta$.

By the definition of $X_p$, we can find a point
$x\in W^{cu}(p)\cap \Omega^c(f)$ such that $d(x,q)<\delta$.
By Lemma~\ref{L0ThmE}, we can take $z\in W^u(x)\cap W^s(q)$.
Since $x\in \Omega^c(f)$,  Theorem D implies that there are
infinitely many points $p_n\in P^c(f)$ such that $p_n\to x$ as $n\to \infty$.
Hence, $d(p_n, q)<\delta$ for all $n$ large enough.
By Lemma~\ref{L0ThmE}, there exist $q_n\in W^c(q)$ and $z_n\in M$ such that
$z_n\in W^u_\ve(p_n)\cap W^s_\ve(q_n)$.
By Lemma~\ref{CloudLemma}, $z_n\in \Omega^c(f)$.
Note that by continuity if $p_n\to x$, then $q_n\to q$ and $z_n\to z$.
We get $z\in \Omega^c(f)$.
Since $z\in W^u(x)$ and  $x\in W^{cu}(p)$, we get $z\in X_p$ by the definition
of $X_p$.
Further, since $z\in W^s(q)$ and $q\in P^c(f)$,
$\{f^n(z)\}_{n\geq 0}$ has at least one accumulation point in $W^c(q)$.
This implies that $W^c(q)\subset X_p$ and we complete the proof.
\end{proof}

\begin{lemma}\label{L2ThmE}
Let $p, q\in P^c(f)$ and $X_p\cap X_q\neq \emptyset$.  Then $X_p=X_q$.
\end{lemma}

\begin{proof}
Since  $X_p\cap X_q\neq \emptyset$, there are points
$x\in W^{cu}(p)\cap \Omega^c(f)$ and $q'\in X_q\cap P^c(f)$
such that $d(x, q')<\delta$.
By Lemma~\ref{L0ThmE}, there exists a point $z\in\Omega^c(f)$
such that $z\in W^u_\ve(x)\cap W^s_\ve(q'_1)$ for some $q'_1\in W^c(q')$.
Let  $n$ be the period of $W^c(q')$.  Then
$$\lim\limits_{i\to+\infty}d\bigl(f^{in}(z), W^c(q')\bigr)=0.
$$
By Lemma~\ref{L1ThmE}, $f^{in}(z)\in X_{q_1}$ for $i$ large enough
and hence $z\in X_q$.

At the same time, we have
$$\lim\limits_{i\to+\infty}d\bigl(f^{-in}(z), W^c(p)\big)=0.
$$
So, $W^c(p)\subset X_q$.

For any $y\in W^{cu}(p)\cap \Omega^c(f)$, one has
$$\lim\limits_{i\to+\infty}d\bigl(f^{-im}(y), W^c(p)\big)=0,
$$
where $m$ is the period of $W^c(p)$.
So $f^{-im}(y)\in B_{\delta,\Omega^c(f)}(W^c(p))\subset B_{\delta,\Omega^c(f)}(X_q)$
for $i$ large enough and hence $y\in X_q$.
Noting that $X_p=\overline{W^{cu}(p)\cap \Omega^c(f)}$,
we have $X_p\subset X_q.$

Similarly, one can get $X_q\subset X_p.$ This completes the proof.
\end{proof}

\hspace*{-5mm}\emph{Acknowledgements}. The second author is supported by NSFC (No: 11471056). The third author is the corresponding author, and is supported by NSFC(No: 11371120), NCET(No: 11-0935), NSFHB(No: A2014205154, A2013205148) and the Plan of Prominent Personnel Selection and Training for the Higher Education Disciplines in Hebei Province (No: BR2-219).

\end{document}